\documentclass[11pt,a4paper]{amsart}

\usepackage[T1]{fontenc}
\usepackage[utf8]{inputenc}
\usepackage[english]{babel}

\usepackage[
    a4paper,
    left=32mm,
    right=32mm,
    top=30mm,
    bottom=30mm
]{geometry}

\usepackage{amsmath}
\usepackage{amssymb}
\usepackage{mathtools}

\usepackage{microtype}
\usepackage{graphicx}
\usepackage{booktabs}
\usepackage{enumitem}

\usepackage[
    colorlinks=false,
    pdfborder={0 0 0.6},
    linkbordercolor={0.65 0.15 0.18},
    citebordercolor={0.10 0.45 0.25},
    urlbordercolor={0.10 0.30 0.65}
]{hyperref}

\usepackage[nameinlink,noabbrev]{cleveref}

\theoremstyle{plain}

\newtheorem{theorem}{Theorem}[section]
\newtheorem{lemma}[theorem]{Lemma}

\theoremstyle{definition}

\newtheorem{definition}[theorem]{Definition}

\theoremstyle{remark}

\newcommand{\R}{\mathbb{R}}

\newcommand{\ddd}{\partial}

\newcommand{\tran}{\mathsf{T}}

\title[Rigidity for the four gradient problem]{A simplified proof of rigidity for the four gradient problem}

\author{Igor Buchowiec}

\address{Faculty of Mathematics, Informatics and Mechanics, University
of Warsaw, ul. Stefana Banacha 2, 02-097 Warsaw, Poland}

\email{i.buchowiec@uw.edu.pl}

\date{November 2026}

\begin{document}

\begin{abstract}
We provide a modified proof of rigidity for the four gradient problem due to Chlebík and Kirchheim (2002). It concerns an inclusion $Du(x)\in\mathcal{K}$ for almost every $x\in\Omega$, a Lipschitz function $u$, a bounded domain $\Omega\subset\R^n$ and a four-element set $\mathcal{K}\subset\R^{m\times n}$ without rank-one connections for $m,n>1$. Rigidity means that every solution $u$ satisfies $Du=\mathrm{const}$ a.e. 
The original proof relies on the regularity theory of solutions to elliptic-type Monge--Amp\`ere equations. Our proof omits the argument based on the regularity theory of Monge--Amp\`ere equation and instead uses only Šverák's convexity theorem.

\end{abstract}



\maketitle

\section{Introduction}
\label{sec:introduction}
We consider the following differential inclusion for a bounded Euclidean domain $\Omega\subset\R^n$, Lipschitz mapping $u\colon\Omega\to\R^m$, subset $\mathcal{K}$ of real $m\times n$ matrices, and integers $m,n>1$:
\begin{equation}\label{zag}
Du(x)\in \mathcal{K}\text{ for a.e. }x\in\Omega.
\end{equation}
We say that two distinct matrices $A,B$ are rank-one connected if $\text{rank}(A-B)=1$.

The $N$-gradient problem asks whether the absence of rank-one connections in a set $\mathcal{K}\subset\R^{m\times n}$ with $N$ matrices forces every Lipschitz solution of \eqref{zag} to be affine. The $N$-gradient problem is called \emph{rigid} if for every set $\mathcal{K}$ with $N$ elements each Lipschitz solution to \eqref{zag} is affine.

Rigidity for the two gradient problem was established by Ball and James \cite{balljames}. Proofs of rigidity for three gradients can be found in the works of Müller \cite{muller} and Kirchheim with Chlebík \cite{chlebik}.
The case of four gradients required further arguments and was resolved by Chlebík and Kirchheim \cite{chlebik} in 2002. They proved that four gradient problem is rigid (see Theorem \ref{czterygrad}).

The restriction to four gradients is essential. The analogous problem with five matrices is nonrigid in general. Kirchheim and Preiss founded a set $\mathcal{K}$ with five symmetric matrices and no rank-one connections for which there exists a nonaffine Lipschitz mapping $u$ satisfying $Du(x)\in\mathcal{K}$ almost everywhere. Their exemplary set $\mathcal{K}$ is given explicitly in \cite{kirch}. The nonaffine solution is constructed using convex integration methods.

The absence of rank-one connections is crucial and makes this result nontrivial. Indeed, following \cite{muller}, let $a\in\R^m$ and $b\in\R^n$ be nonzero. Define the $m\times n$ matrix $a\otimes b$ by $(a\otimes b)_{ij}=a_ib_j$. Suppose that there are $A,B\in\mathcal{K}$ satisfying $\mathrm{rank}(B-A)=1$. Hence, we can write $B-A=a\otimes b$. Then, set
\[
u(x)=Ax+ah(\langle x,b\rangle)
\]
and choose $h\colon\R\to\R$ to be a Lipschitz function satisfying $h'(t)\in\{0,1\}$ almost everywhere. It follows that $Du(x)\in\{A,B\}$ for almost every $x\in\Omega$. Taking $h$ as a piecewise linear map with both derivatives occurring on subsets of positive measure gives a nonaffine solution. This shows that the presence of a rank-one connection reduces the construction of nontrivial solutions to a one-dimensional problem.

The original proof of the four-gradient theorem relies on three highly nontrivial ingredients: openness  and positivity of the Jacobian  of a nonconstant quasiregular mapping (Theorem \ref{jakobianqr}), convexity of planar solutions of the Monge--Amp\`ere equation (Theorem \ref{wypuklosc}) and regularity theory of such solutions (Theorem \ref{caffa}).

In this paper, we revisit the original proof and simplify its final part concerning regularity theory for Monge--Amp\`ere equation. After reducing the problem to the planar case, we make the algebraic reduction explicit by giving formulas for the transformations leading to symmetric matrices with the same determinant. When this determinant is positive, Šverák's convexity theorem (\cite{sverak}, see Theorem \ref{wypuklosc} below) applies. It turns out that convexity is sufficient to obtain rigidity of the problem.

The paper is organized as follows. After stating the main result, we recall the results concerning reductions and rigidity of gradient problems with two and three gradients, together with the required facts about quasiregular mappings. Then, we proceed with convexity and regularity results for the Monge--Amp\`ere equation. In particular, we prove an elementary lemma on convexity that provides the final step of our argument. In the last section of the paper, we present the simplified proof in detail.

\section{Auxiliary tools and rigidity results}

Throughout, $\Omega$ denotes a bounded Euclidean domain, and almost everywhere statements refer to Lebesgue measure. For Lipschitz mappings, classical and weak derivatives are identified almost everywhere. For a scalar-valued function $u$ with Lipschitz continuous gradient, $D^2u$ denotes the almost everywhere defined derivative of $\nabla u$, which agrees almost everywhere with the weak Hessian and is symmetric.

Our aim is to present a simplified proof of the following rigidity theorem originally proved by Chlebík and Kirchheim.

\begin{theorem}\label{czterygrad}\textup{(\cite{chlebik})}
    Let $\mathcal{K}=\{A_1,A_2,A_3,A_4\}$ be a subset of $\R^{m\times n}$ with no rank-one connections, and let $\Omega\subset\R^n$ be a bounded domain. If $u\colon\Omega\to\R^m$ is Lipschitz and $Du(x)\in\mathcal{K}$ for almost every $x\in\Omega$, then $Du(x)=\mathrm{const}$ almost everywhere.
\end{theorem}

As discussed before, the absence of rank-one connections is essential, and the conclusion fails in general for sets containing five matrices. Our proof retains the quasiregular and convexity arguments of the original proof, while replacing its final use of elliptic regularity with a direct sign contradiction.

\subsection{Auxiliary results}

The following lemma, due to Chlebík and Kirchheim, states that rigidity of the finite gradient problem in a planar bounded domain case implies rigidity in higher dimensions.

\begin{lemma}\textup{(\cite{kirch}, \cite{chlebik})}\label{lemacik}
    Let $\Tilde{\mathcal{K}}\subset \R^{m\times n}$ be a finite subset with no rank-one connections. Let $\Omega\subset\R^n$ be a bounded domain. Suppose that there exists a non-affine Lipschitz function $\Tilde{u}\,\colon\Omega\to\R^m$ such that $D\Tilde{u}(x)\in\Tilde{\mathcal{K}}$ for almost every $x\in\Omega$. Then there exist linear maps $g\,\colon\R^2\to\R^n$ and $h\,\colon\R^m\to\R^2$, identified respectively with their matrices $G\in\R^{n\times 2}$ and $H\in \R^{2\times m}$ in the corresponding standard bases, and a Lipschitz function $u\,\colon[0,1]^2\to\R^2$ satisfying:
    \begin{enumerate}[label = (\roman*)]
        \item $\mathcal{K}=\{HAG\colon\,A\in\Tilde{\mathcal{K}}\}
        \subset \R^{2\times 2}$ has no rank-one connections;
        \item $Du(x)\in \mathcal{K}$ for almost every $x\in[0,1]^2$;
        \item $u$ is not affine.
    \end{enumerate}
\end{lemma}

We recall rigidity results for two and three gradients.
Their proofs are omitted.

\begin{theorem}\label{dwagrad}
    Let $\mathcal{K}$ be a subset of $\R^{m\times n}$ with at most three elements and no rank-one connections. If $u\colon\Omega\subset\R^n\to\R^m$ is Lipschitz and $Du(x)\in \mathcal{K}$ for almost every $x\in\Omega$,
    then $Du(x)=\mathrm{const}$ almost everywhere.
\end{theorem}

The proof of rigidity for the two gradient problem is due to Ball and James (see \cite{balljames}). The proof of rigidity with $\mathcal{K}$ consisting of three matrices can be found in \cite{muller} and \cite{kirch}.


\subsection{Quasiregular mappings and applications to gradient problems}
We start with the definition of a quasiregular mapping in dimension two. We restrict ourselves to Lipschitz quasiregular mappings since this is the case in gradient problems. 

\begin{definition}
    The Lipschitz continuous mapping $f\colon\Omega\to\R^2$, where $\Omega\subset\R^2$, is called \emph{quasiregular} if it satisfies
    \[
    \|Df(x)\|^2\leq K\det Df(x)
    \]
    for almost every $x\in\Omega$ and some constant $K\geq 1$.
\end{definition}

The following classical result due to Reshetnyak \cite{reshet} is helpful in gradient problems.

\begin{theorem}\label{jakobianqr}\textup{(\cite{reshet})}
    If $f:\Omega\to\R^2$ is a nonconstant quasiregular
    mapping, then $\det Df(x)>0$ for almost every
    $x\in\Omega$. Moreover, each nonconstant quasiregular mapping is an open and discrete mapping.
\end{theorem}

The first part of this theorem is used in the gradient problems as a tool for reducing the number of gradients. Indeed, let us consider the $N$-gradient problem with one gradient being the zero matrix and the remaining matrices having positive determinants. Then the solution is quasiregular. Therefore, we conclude that the solution is either constant or attains nontrivially at most ($N-1$)-gradients. For our purposes, the case $N=4$ will be sufficient.

\subsection{Convexity and regularity for the Monge--Amp\`ere equation}

We collect the convexity and regularity results used
in the original proof of the four gradient theorem. In the last section we will present the proof without regularity theory. We include references to complete proofs of these results.

The first result, concerning convexity, was proved by Šverák in his unpublished manuscript \cite{sverak}. The proof can be found in \cite{lewicka}. 

\begin{theorem}\textnormal{(Theorem 2.3 in \cite{lewicka}; \cite{sverak})}\label{wypuklosc}
    Let $\Omega$ be an open, convex and bounded subset of $\R^2$. Assume that $f\colon\Omega\to\R$ is strictly positive, meaning that $f(x)\geq c_0>0$ for some constant $c_0$ and almost every $x\in\Omega$. Suppose that $u\colon\Omega\to\R^2$ has Lipschitz gradient and satisfies the Monge--Amp\`ere equation
    \begin{equation}\label{monge}
        \det D^2u(x)=f(x)
    \end{equation}
    for almost every $x\in\Omega$. Then, possibly after a global sign change of $u$ ($u\mapsto -u$), the function $u$ is convex on $\Omega$.
\end{theorem}

The following two results will not be used in our modified proof, although they are presented for completeness. 

Firstly, in the planar case we obtain not only convexity, but also strict convexity.
\begin{theorem}\textnormal{(Remark 2.4 in \cite{mooney})}\label{scisla}
    Let $\Omega\subset\R^2$ be open, bounded and convex. If a convex function $u\colon\Omega\to\R$ with a Lipschitz continuous gradient satisfies the \emph{Monge--Amp\`ere-type} inequality
    \[
        \det D^2u(x)\geq 1
    \]
    for almost every $x\in\Omega$, then $u$ is strictly convex.
\end{theorem}

Strict convexity allows to use the regularity theory for the elliptic Monge--Amp\`ere equation and obtain smoothness of the solution if the right-hand side of \eqref{monge} is smooth. 

The following result by Caffarelli \cite{caffa} was used in the original proof of Theorem \ref{czterygrad}.

\begin{theorem}\textnormal{(\cite{caffa})}\label{caffa}
    If the function $f$ specified in \textnormal{Theorem \ref{wypuklosc}} belongs to $\mathcal{C}^\infty(\Omega)$, then every strictly convex $u\colon\Omega\to\R^2$ with Lipschitz gradient solving the equation \textnormal{(\ref{monge})} also belongs to $\mathcal{C}^\infty(\Omega)$.
\end{theorem}

Instead of the regularity theory, the following elementary consequence of convexity will suffice in our proof of rigidity for the four gradient problem.

\begin{lemma}\label{minory}
    Let $\Omega\subset\R^2$ be an open ball. If a function $u\colon\Omega\to\R$ with a Lipschitz continuous gradient is convex (concave, respectively), then $D^2u(x)$ is a positive (negative) semidefinite matrix for almost every $x\in\Omega$. In particular, the upper-left entry of $D^2u(x)$ is nonnegative (nonpositive)
    for almost every $x\in\Omega$.
\end{lemma}

\begin{proof}
    We consider the convex case. The concave one follows by taking $-u$ instead of $u$.
    
    The second part of the lemma follows from Sylvester's criterion for semidefiniteness of a symmetric matrix. It remains to prove the first part. If $u\in\mathcal{C}^2(\Omega)$, the result is immediate from one of the equivalent definitions of convexity. Without this additional regularity, for $\varepsilon>0$ we consider the convolution $u_{\varepsilon}(x)\coloneqq(u\ast\phi_{\varepsilon})(x)$ with a standard mollifier $\phi_{\varepsilon}$.

    This convolution is well defined on $\Omega_{\delta}$, the set of points in $\Omega$ whose distance from $\ddd\Omega$ is greater than $\delta>0$. From now on, $\delta>0$ is arbitrary and $\varepsilon<\delta$. Since $u$ is locally integrable (being continuous, as it has a Lipschitz continuous gradient), the convolution is well defined on $\Omega_{\delta}$, and $u_{\varepsilon}(x)\to u(x)$ as $\varepsilon\to 0^+$ for almost every $x\in\Omega_{\delta}$.

    Since $u$ has a Lipschitz continuous gradient, we also obtain pointwise convergence $D^2u_{\varepsilon}(x)\to D^2u(x)$ as $\varepsilon\to 0^+$ for almost every $x\in\Omega_{\delta}$ by the identity $D^2u_{\varepsilon}=(D^2u)\ast\phi_{\varepsilon}$ and $D^2u\in L^1_{\text{loc}}(\Omega_{\delta})$. Moreover, a straightforward calculation shows that $u_{\varepsilon}$ is smooth and convex. By the characterization of convexity for functions of class $\mathcal{C}^2(\Omega_{\delta})$, for every $y=[y_1,y_2]^{\tran}\in\R^2$ and every $x\in\Omega_{\delta}$ we have
    \[
        y^{\tran}D^2u_{\varepsilon}(x)y\geq 0.
    \]
    The almost everywhere pointwise convergence of the Hessians yields
    \[
        y^{\tran}D^2u(x)y\geq 0
    \]
    for almost every $x\in\Omega_{\delta}$ and every $y\in\R^2$. Thus, $D^2u(x)$ is positive semidefinite  for almost every $x\in\Omega_{\delta}$. Applying this to $\delta=1/n$ and using $\Omega=\bigcup_{n=1}^{\infty}\Omega_{1/n}$, we obtain the assertion for almost every $x\in\Omega$. This completes the proof of the lemma.
\end{proof}

\section{A simplified proof of the main result}
\label{sec:main-results}

The following theorem by Chlebík and Kirchheim is a crucial part of the original proof of rigidity for the four gradient problem. It will be used in our simplified proof.

\begin{theorem}\textnormal{(\cite{chlebik},\cite{kirch})}\label{kirch}
    Let $\mathcal{K}$ be a four-element subset of $\R^{2\times 2}_{\text{sym}}\cap \{A\colon\det A=-d\}$ with no rank-one connections, where $d\geq 0$. If $u\colon\Omega\to\R^2$ is Lipschitz and $Du(x)\in \mathcal{K}$ for almost every $x\in\Omega$, then $u$ is affine. In particular, $Du(x)=\mathrm{const}$ almost everywhere.
\end{theorem}

Its proof uses some rather involved techniques and relies on openness of nonconstant quasiregular mappings (Theorem \ref{jakobianqr}) and geometric measure theory. Notably, this result does not rely on either convexity results or regularity theory for the Monge--Amp\`ere equation.

\begin{proof}[Proof of Theorem \ref{czterygrad}]
The reasoning is divided into five steps.
\begin{center}
    \textsc{Step 0. (reduction)}
\end{center}
By Lemma~\ref{lemacik}, it suffices to consider the case $m=n=2$. We may assume that $\Omega$ is an open ball, $\Omega=B(x_{\omega},\omega)$. Indeed, if $u$ is affine on every ball contained in $\Omega$, the affine representations agree on overlaps, and a chain of overlapping balls in the bounded domain $\Omega$ shows that $u$ is affine on the whole domain. Translating and rescaling by setting $v(x)\coloneqq 1/\omega\cdot u(x_{\omega}+\omega x)$ we preserve the gradient inclusion and obtain a problem on the unit ball. Retaining the notation $u$ for the rescaled function, we henceforth take $\Omega=B(0,1)$.

The proof proceeds by contradiction. Therefore, we assume that $u$ is not affine. Then by Theorem \ref{dwagrad}, each of the four gradients is attained on a set with positive measure in $\Omega$.

\begin{center}
    \textsc{Step 1. (analysis based on quasiregularity)}
\end{center}
By subtracting a linear mapping and applying an invertible outer linear transformation the problem \eqref{zag} can be reduced to
\[
Du(x)\in\{0,I,B_1,B_2\}
\]
for almost every $x\in\Omega$. Invertible linear transformations preserve the absence of rank-one connections since the product of matrices with full rank has full rank.

Thus, $\text{rank}(B_1)=\text{rank}(B_2)=\text{rank}(B_1-I)=\text{rank}(B_2-I)=\text{rank}(B_2-B_1)=2$, since $\mathcal{K}$ had no rank-one connections. In particular, $\det(B_1),\det(B_2)\neq 0$. There are four possible cases:
\begin{enumerate}
    \item $\det(B_1),\det(B_2)>0$;
    \item $\det(B_1),\det(B_2)<0$;
    \item $\det(B_1)>0$, $\det(B_2)<0$;
    \item $\det(B_1)<0$, $\det(B_2)>0$.
\end{enumerate}
In the first case, $u$ must be quasiregular. Thus, by Theorem \ref{jakobianqr}, $u$ is either constant or $u$ uses at most three gradients. Theorem~\ref{dwagrad} then implies that $u$ is affine.

The third case is equivalent to the fourth one after relabelling the matrices $B_i$ and the second case reduces to the fourth one. Indeed, multiplying on the left by $B_1^{-1}$ (an outer transformation), we obtain the set $\{0,I,B_1^{-1},B_1^{-1}B_2\}$. The determinants have the same sign pattern as in the fourth case.

From now on, without loss of generality, we consider only the last case, that is, the case of the zero matrix, the identity matrix and two matrices with determinants of negative and positive sign denoted by $B_1$, $B_2$ respectively.

We use the fact that a $2\times 2$ matrix with negative determinant is diagonalizable over $\R$. We diagonalize $B_1$ by applying inner and outer transformations $B_1\mapsto Q^{-1}B_1Q$ for some invertible matrix $Q=(q_{ij})_{i,j=1,2}$. The precise transformation on $u$ takes the form $u(x_1,x_2)\mapsto Q^{-1}u\left(Q(x_1,x_2)\right)=Q^{-1}u(q_{11}x_1+q_{12}x_2,q_{21}x_1+q_{22}x_2)$.

Since both $0$ and $I$ remain unchanged under the diagonalization of $B_1$, after this transformation our problem can be written as follows:
\[
Du(x)\in \mathcal{K}\coloneqq\left\{0,I,B\coloneqq\begin{bmatrix}
    \alpha & 0\\
    0 & \beta
\end{bmatrix},C\coloneqq Q^{-1}B_2 Q\right\}
\]
for almost every $x\in\Omega$, where $\alpha\beta<0$, $\text{rank}(B-I)=\text{rank}(C)=\text{rank}(C-I)=\text{rank}(C-B)=2$ and $\det(C)>0$. The domain remains bounded, open and convex under this change of variables. Let \( C=\bigl[\begin{smallmatrix}
a & b\\
c & d
\end{smallmatrix}\bigr] \).
Two cases remain
\begin{equation*}
    \text{a) $b=0$ or $c=0$},
    \quad
    \text{b) $b\neq 0$ and $c\neq 0$.}
\end{equation*}
The (triangular) case a) will be treated in \textsc{Step 2} of the proof. In the last two steps of the proof we will resolve the b) case.
\begin{center}
    \textsc{Step 2. (case \textnormal{a}) - triangular matrices)}
\end{center}
Suppose that $c=0$. Set
\[
J\coloneqq\{x_2\in\R:(x_1,x_2)\in\Omega \text{ for some }x_1\in\R\},
\]
and, for $x_2\in J$, define
\[
I_{x_2}\coloneqq\{x_1\in\R:(x_1,x_2)\in\Omega\},
\quad
L_{x_2}\coloneqq I_{x_2}\times\{x_2\}.
\]
Since $\Omega$ is open and convex, each $I_{x_2}$ and $J$ are open intervals.

Notice that $\ddd u_2/\ddd x_1=0$ almost everywhere. By absolute continuity along horizontal lines (\cite{mazya}), $u_2$ is constant on $L_{x_2}$ for almost every $x_2\in J$. Continuity then gives the same conclusion for every $x_2\in J$. Thus $u_2$ is independent of $x_1$, and we may write $u_2(x_1,x_2)=h(x_2)$ in $\Omega$. Let
\[
F\coloneqq
\{x=(x_1,x_2)\in\Omega: u\text{ is differentiable at }x \text{ and }Du(x)\in \mathcal{K}\}.
\]
This is a set of full measure in $\Omega$. We first show that $Du$ is independent of $x_1$ on the intersection of each horizontal line with $F$.

Fix two arbitrary points $(x^0_1,x^0_2),(\tilde{x}^0_1,x^0_2)\in F$ on the same horizontal line, and suppose that the assertion fails. We can then write
\[
Du(x^0_1,x^0_2)
=M_1\coloneqq\begin{bmatrix}
    a_1 & b_1\\
    0 & d_1
\end{bmatrix},
\quad
Du(\tilde{x}^0_1,x^0_2)
=M_2\coloneqq\begin{bmatrix}
    a_2 & b_2\\
    0 & d_2
\end{bmatrix},
\]
where $M_1,M_2\in \mathcal{K}$, $M_1\neq M_2$, $a_1,a_2\in\{0,1,\alpha,a\}$, $b_1,b_2\in\{0,b\}$ and $d_1,d_2\in\{0,1,\beta,d\}$. However, $u_2(x_1,x_2)=h(x_2)$, so
\[
d_1
=\frac{\ddd u_2}{\ddd x_2}(x^0_1,x^0_2)
=h'(x^0_2)
=\frac{\ddd u_2}{\ddd x_2}(\tilde{x}^0_1,x^0_2)
=d_2.
\]
The derivative $h'(x^0_2)$ exists because $u$ is differentiable at the selected points. This is impossible: it would give $\text{rank}(M_1-M_2)<2$, contrary to the assumptions on rank-one connections in $\mathcal{K}$.

By Fubini's theorem, for almost every $x_2\in J$, the section $F\cap L_{x_2}$ has full one-dimensional measure in $L_{x_2}$. On each such section, the derivative of the mapping $x\mapsto u(x_1,x_2)$ is constant almost everywhere. Since this function is Lipschitz, absolute continuity implies that it is affine on the whole interval $I_{x_2}$. We claim that $u$ is affine on every horizontal section $L_{x_2}$.

Fix an arbitrary $x_2^0\in J$. Since a set of full measure in an interval is dense, there exists a sequence $x_2^n\to x_2^0$ such that $u$ is affine on every $L_{x_2^n}$. Thus for every $x_1^n\in I_{x_2^n}$ we have
\[
u(x_1^n,x_2^n)=A_nx_1^n+B_n
\]
for some $A_n,B_n\in\R^2$.

Choose two distinct points $x_1,\tilde{x}_1\in I_{x_2^0}$. Since $\Omega$ is open, both $(x_1,x_2^n)$ and $(\tilde{x}_1,x_2^n)$ belong to $\Omega$ for all sufficiently large $n$. By Lipschitz continuity of $u$ and the identities
\[
A_n=\frac{u(x_1,x_2^n)-u(\tilde{x}_1,x_2^n)}{x_1-\tilde{x}_1},
\quad
B_n=u(x_1,x_2^n)-A_nx_1,
\]
the sequences $(A_n)$ and $(B_n)$ are bounded. Passing to a common subsequence, denoted in the same way, we may assume that $A_n\to A_0$ and $B_n\to B_0$.

For every fixed $\hat{x}_1\in I_{x_2^0}$, openness of $\Omega$ also ensures that $(\hat x_1,x_2^n)\in\Omega$ for all sufficiently large $n$. Consequently,
\[
u(\hat x_1,x_2^0)
=\lim_{n\to\infty}u(\hat x_1,x_2^n)
=\lim_{n\to\infty}(A_n\hat x_1+B_n)
=A_0\hat x_1+B_0.
\]
Since $x_2^0$ was arbitrary, this proves the claimed property that $u$ is affine on every horizontal section $L_{x_2}$.

We can therefore write
\[
u(x_1,x_2)=S(x_2)x_1+P(x_2)
\]
for some functions $S,P:J\to\R^2$ and all $(x_1,x_2)\in\Omega$. We show that $S$ is continuous.

Fix $x_2\in J$ and a sequence $(x_2^n)$ in $J$ converging to $x_2$. Choose two distinct points $x_1,\tilde x_1\in I_{x_2}$. By openness of $\Omega$, the points $(x_1,x_2^n)$ and $(\tilde x_1,x_2^n)$ belong to $\Omega$ for all sufficiently large $n$. By the continuity of $u$,
\[
S(x_2^n)(x_1-\tilde x_1)=u(x_1,x_2^n)-u(\tilde x_1,x_2^n)\longrightarrow u(x_1,x_2)-u(\tilde x_1,x_2)
=S(x_2)(x_1-\tilde x_1).
\]
Since $x_1-\tilde x_1\neq 0$, we obtain $S(x_2^n)\to S(x_2)$. Thus $S$ is continuous.

By the definition of the directional derivative, for almost every $(x_1,x_2)\in\Omega$ we have
\[
S(x_2)=Du(x_1,x_2)\cdot e_1
\in
\left\{
\begin{bmatrix}
    0\\
    0
\end{bmatrix},
\begin{bmatrix}
    1\\
    0
\end{bmatrix},
\begin{bmatrix}
    \alpha\\
    0
\end{bmatrix},
\begin{bmatrix}
    a\\
    0
\end{bmatrix}
\right\}
\eqqcolon \mathcal{M}.
\]

It remains to show that $S$ is constant. 

By Fubini's theorem, $S(x_2)\in \mathcal{M}$ for almost every $x_2\in J$. Since $S$ is continuous and $M$ is closed, this inclusion holds for every $x_2\in J$. The interval $J$ is connected and $\mathcal{M}$ is finite, so $S$ must be constant everywhere.

The absence of rank-one connections implies that the numbers $0,1,\alpha,a$ are pairwise distinct. Indeed, the difference of any two distinct matrices in $\mathcal{K}$ is upper triangular and has rank two, so both its diagonal entries are nonzero. Thus the first column uniquely determines a matrix in $\mathcal{K}$. Since $Du(x_1,x_2)\cdot e_1=S$ is constant, it follows that $Du$ is constant almost everywhere in $\Omega$. As $\Omega$ is connected, $u$ is affine. 

The case $b=0$ is treated in the same way,
with horizontal lines replaced by vertical ones
and the roles of the two components of $u$ interchanged.\footnote{Although rigidity for triangular matrices is known (see \cite{muller}, \cite{chlebik}), we present a detailed, self-contained proof for clarity.}
    
We are left with the case b). We treat it in the next two steps.
\begin{center}
\textsc{Step 3. (case \textnormal{b}) - reduction to symmetry and the same determinant)}
\end{center}
Our aim is to reduce the case b) to the situation as in the Monge--Amp\`ere equation. First, we show that the analysis can be reduced to symmetric matrices. Multiplying on the left by the invertible matrix \( \bigl[\begin{smallmatrix}
1 & 0\\
0 & b/c
\end{smallmatrix}\bigr] \) we obtain the inclusion
\[Du(x)\in\left\{0,\begin{bmatrix}
    1 & 0\\
    0 & \frac{b}{c}
\end{bmatrix},\begin{bmatrix}
    \alpha & 0\\
    0 & \frac{b\beta}{c}
\end{bmatrix},\begin{bmatrix}
    a & b\\
    b & \frac{bd}{c}
\end{bmatrix}\right\}\]
for almost every $x\in\Omega$. By our assumptions we have $\alpha\neq 1$, $t\coloneqq\frac{b}{c}\neq 0$ and $\alpha\beta<0$. In particular, $\alpha\neq\beta$. It remains to make the determinants equal by adding a common symmetric matrix. It turns out that for some $s\in\R$ the system
\begin{equation}\label{wyznacznik}
    \det\begin{bmatrix}
    p & q\\
    q & r
\end{bmatrix}=\det\begin{bmatrix}
    p+1 & q\\
    q & r+t
\end{bmatrix}=\det\begin{bmatrix}
    p+\alpha & q\\
    q & r+\beta t
\end{bmatrix}=\det\begin{bmatrix}
    p+a & q+y\\
    q+y & r+td
\end{bmatrix}=s
\end{equation}
has a real solution for $p,q,r$ under those constraints.
The solvability of \eqref{wyznacznik} is equivalent to the solvability of the following system with variables $p,q,r,s\in\R$ and constant $\alpha,\beta,t\in\R$.
\begin{equation}\label{uklad}
\left\{
\begin{array}{l}
pr-q^2=s,\\
(p+1)(r+t)-q^2=s,\\
(p+\alpha)(r+\beta t)-q^2=s,\\
(p+a)(r+td)-(b+q)^2=s.
\end{array}
\right.
\end{equation}
We show how to solve this system using our assumptions and then use precise formulas to conclude in \textsc{Step 4}.
By subtracting the first equation from the second one we get $r=-t(p+1)$. Next, subtracting the first equation from the third one and putting $r=-t(p+1)$ we obtain a linear equation in $p$ with a solution $p=\frac{\alpha-\alpha\beta}{\beta-\alpha}$. It is well-defined by our assumptions on $\alpha$ and $\beta$. Taking $p$ and $q$ in those forms into the first equation yields the value of $s$, where $s$ depends on $t,\alpha,\beta$ and on $q$. It remains to find $q$. Subtracting the first equation from the last one, reduces to the following equation for $q$ with $p$ and $r$ independent on $q$
\[ptd+ar+atd-b^2-2bq=0.\]
It is a linear equation in $q$, since $b\neq 0$. This shows that the system \eqref{uklad} has a real solution. In the final step of the proof we conclude by analyzing the sign of the common determinant.
\begin{center}
\textsc{Step 4. (case \textnormal{b}) - contradiction with convexity)}
\end{center}
By \textsc{Step 3} we know that it remains to consider the inclusion $Du(x)\in \mathcal{L}$ for almost every $x\in\Omega$ and for the four-element subset $\mathcal{L}\subset \R^{2\times 2}$ consisting of symmetric matrices with the same determinant $s$.
    
If $s\leq 0$, then by Theorem \ref{kirch} $u$ is affine and the proof is finished.
    
If $s>0$, then thanks to symmetry and convexity of $\Omega$ we are left with the Monge--Amp\`ere equation, i.e.,  $\det D^2v(x)=s$ for scalar-valued function $v\colon\Omega\to\R$ with Lipschitz continuous gradient ($D^2 v(x)=Du(x)$). By Theorem \ref{wypuklosc} the function $v$ is either convex or concave on $\Omega$. This result applies since $\Omega$ is an open ball, thus a convex set. In both cases, Lemma \ref{minory} enforces the product of two upper-left cells $(p+\alpha)$, $(p+1)$ to be nonnegative since every gradient is attained on a set of positive measure.
However, as $p=\frac{\alpha(1-\beta)}{\beta-\alpha}$
\[(p+\alpha)(p+1)=\frac{\alpha(1-\alpha)}{\beta-\alpha}\cdot\frac{\beta(1-\alpha)}{\beta-\alpha}=\alpha\beta\left(\frac{1-\alpha}{\beta-\alpha}\right)^2<0.
\] 
We used the fact that $\alpha\beta<0$, $\alpha\neq\beta$ and $\alpha\neq 1$. This contradiction proves the result.
\end{proof}

\section*{Acknowledgements}
The author would like to thank Agnieszka Kałamajska for valuable discussions concerning this study. This research is co-financed by the National Science Centre, Poland grant number 2023/51/D/ST1/02907.

\section*{Disclaimer on the use of LLMs.}
This paper is entirely the author’s own work, written independently in 2025 without any assistance from Large Language Models. The mathematical arguments and original text were developed entirely by the author. No AI-generated text, proofs, or mathematical arguments have been incorporated into the manuscript.

\end{document}